\documentclass[12pt,letterpaper,reqno]{amsart}

\newcount\timehh\newcount\timemm

\timehh=\time

\divide\timehh by 60 \timemm=\time

\usepackage{times}

\usepackage[T1]{fontenc}

\usepackage{mathrsfs}

\usepackage{latexsym}

\usepackage[dvips]{graphics}

\usepackage{epsfig}

\usepackage{amsmath,amsfonts,amsthm,amssymb,amscd}

\input amssym.def

\input amssym.tex

\newcommand\be{\begin{equation}}

\newcommand\ee{\end{equation}}

\newcommand\bea{\begin{eqnarray}}

\newcommand\eea{\end{eqnarray}}

\newcommand\bi{\begin{itemize}}

\newcommand\ei{\end{itemize}}

\newcommand\ben{\begin{enumerate}}

\newcommand\een{\end{enumerate}}

\newcommand\bc{\begin{center}}

\newcommand\ec{\end{center}}

\newcommand\ba{\begin{array}}

\newcommand\ea{\end{array}}

\numberwithin{equation}{section}

\begin{document}

\title{Irrationality from The Book}% (alternate: Rational Irrational Proofs)}

\author{Steven J. Miller}\email{Steven.J.Miller@williams.edu}

\address{Department of Mathematics and Statistics, Williams College, Williamstown, MA 01267}

\author{David Montague}\email{davmont@umich.edu}

%   {
%      \protect \protect\sc\today\ --
%      \ifnum\timehh<10 0\fi\number\timehh\,:\,\ifnum\timemm<10
%0\fi\number\timemm
%      \protect \, \, \protect \bf DRAFT -- DO NOT DISTRIBUTE
%   }

\address{Department of Mathematics, University of Michigan, Ann Arbor, MI 48109}

\date{\today}

\thanks{This paper was inspired by Margaret Tucker's senior colloquium talk at Williams College (February 9, 2009), where she introduced the first named author to Tennenbaum's wonderful proof of the irrationality of $\sqrt{2}$, and a comment during the talk by Frank Morgan, who wondered if the method could be generalized to other numbers. We thank Peter Sarnak for pointing out the reference \cite{Cw}. Some of the work was done at the 2009 SMALL Undergraduate Research Program at Williams College and the 2009 Young Mathematicians Conference at The Ohio State University; it is a pleasure to thank the NSF (Grant DMS0850577), Ohio State and Williams College for their support; the first named author was also supported by NSF Grant DMS0600848.}

\maketitle

A right of passage to theoretical mathematics is often a proof of the irrationality of $\sqrt{2}$, or at least this is where a lot of our engineering friends part ways, shaking their heads that we spend our time on such pursuits. While there are numerous applications of irrational numbers (these occur all the time in dynamical systems), the purpose of this note is not to sing their praise, but to spread the word on a truly remarkable geometric proof and its generalizations.

Erd\"os used to say that G-d kept the most elegant proofs of each mathematical theorem in `The Book', and while one does not have to believe in G-d, a mathematician should believe in The Book. We first review the standard proof of the irrationality of $\sqrt{n}$. We then give the elegant geometric proof of the irrationality of $\sqrt{2}$ by Stanley Tennenbaum (discovered in the early 1950s \cite{Te}, and which first appeared in print in John H. Conway's article in Power \cite{Co}), which is worthy of inclusion in The Book. We then generalize this proof to $\sqrt{3}$ and $\sqrt{5}$, and invite the reader to explore other numbers.

%%%%%%%%%%%%%%%%%%%%%%%%%%%%%%%%%%%%%%%%%%%%%%%%%%%%%%%%%%%%%%%%%%%%%%%%%%%%%%%%%%%%%%%%%%%%%%%%%%%%%%%%%%%%%%%%

%%%%%%%%%%%%%%%%%%%%%%%%%%%%%%%%%%%%%%%%%%%%%%%%%%%%%%%%%%%%%%%%%%%%%%%%%%%%%%%%%%%%%%%%%%%%%%%%%%%%%%%%%%%%%%%%

\section{Standard proof that $\sqrt{2}$ is irrational}

To say $\sqrt{n}$ is irrational means there are no integers $a$ and $b$ such that $\sqrt{n} = a/b$. There are many proofs; see \cite{Bo} for an extensive list. One particularly nice one can be interpreted as an origami construction (see proof 7 of \cite{Bo} and the references there, and pages 183--185 of \cite{CG} for the origami interpretation). Cwikel \cite{Cw} has generalized these origami arguments to yield the irrationality of other numbers as well.

We give the most famous proof in the case $n=2$. Assume there are integers $a$ and $b$ with $\sqrt{2} = a/b$; note we may cancel any common factors of $a$ and $b$ and thus we may assume they are relatively prime. Then $2b^2 = a^2$ and hence $2|a^2$. We want to conclude that $2|a$. If we know unique factorization\footnote{Unique factorization states that every integer can be written as a product of prime powers in a unique way.}, the proof is immediate. If not, assume $a = 2m+1$ is odd. Then $a^2 = 4m^2 + 4m + 1$ is odd as well, and hence not divisible by two.\footnote{Some texts, such as \cite{HW}, state that the Greeks argued along these lines, which is why they stopped their proofs at something like the irrationality of $\sqrt{17}$, as they were looking at special cases and not using unique factorization.} We therefore may write $a = 2r$ with $0 < r < a$. Then $2b^2 = a^2 = 4 r^2$, which when we divide by 2 gives $b^2 = 2r^2$. Arguing as before, we find that $2|b$, so we may write $b=2s$. But now we have $2|a$ and $2|b$, which contradicts $a$ and $b$ being relatively prime. Thus, $\sqrt{2}$ is irrational.

%%%%%%%%%%%%%%%%%%%%%%%%%%%%%%%%%%%%%%%%%%%%%%%%%%%%%%%%%%%%%%%%%%%%%%%%%%%%%%%%%%%%%%%%%%%%%%%%%%%%%%%%%%%%%%%%

%%%%%%%%%%%%%%%%%%%%%%%%%%%%%%%%%%%%%%%%%%%%%%%%%%%%%%%%%%%%%%%%%%%%%%%%%%%%%%%%%%%%%%%%%%%%%%%%%%%%%%%%%%%%%%%%

\section{Tennenbaum's proof}

\begin{figure}[h]
\begin{center}
\scalebox{.7}{\includegraphics{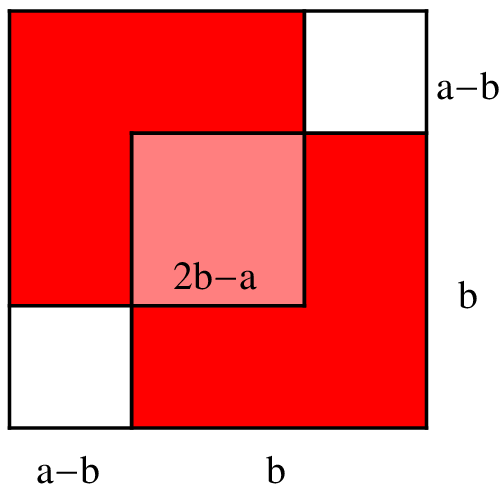}}
\caption{\label{fig:sqrt2} Geometric proof of the irrationality of $\sqrt{2}$.}
\end{center}\end{figure}

We now describe Tennenbaum's wonderful geometric proof of the irrationality of $\sqrt{2}$. Suppose that $(a/b)^2 = 2$ for some integers $a$ and $b$. Without loss of generality, we might as well assume that $a$ and $b$ are the smallest such numbers. So far this looks remarkably similar to the standard proof; here is where we go a different route. Again we find $a^2 = 2b^2$, but now we interpret this geometrically as the area of two squares of side length $b$ equals the area of one square of side length $a$. Thus, if we consider Figure \ref{fig:sqrt2},   we see that the total area covered by the squares with side length $b$ (counting the overlapping, pink region twice) is equal to the area of the larger square with side length $a$. Therefore the pink, doubly counted part, which is a square of side length $2b-a$, has area equal to that of the two white, uncovered squares (each of side length $a-b$) combined. In other words, $(2b-a)^2 = 2(a-b)^2$ or $\sqrt{2} = (2b-a)/(a-b)$. But this is a smaller solution\footnote{While it is clear geometrically that this is a smaller solution, we can also see this algebraically by observing that $2b-a < a$. This follows from $b<a$ and writing $2b-a$ as $b - (a-b) < b < a$. Similarly one can show that $a-b < b$ (if not, $a \ge 2b$ so $a/b \ge 2 > \sqrt{2}$).}, contradiction!

%%%%%%%%%%%%%%%%%%%%%%%%%%%%%%%%%%%%%%%%%%%%%%%%%%%%%%%%%%%%%%%%%%%%%%%%%%%%%%%%%%%%%%%%%%%%%%%%%%%%%%%%%%%%%%%%

%%%%%%%%%%%%%%%%%%%%%%%%%%%%%%%%%%%%%%%%%%%%%%%%%%%%%%%%%%%%%%%%%%%%%%%%%%%%%%%%%%%%%%%%%%%%%%%%%%%%%%%%%%%%%%%%

\section{The square-root of $3$ is irrational}

\begin{figure}[h]
\begin{center}
\scalebox{.7}{\includegraphics{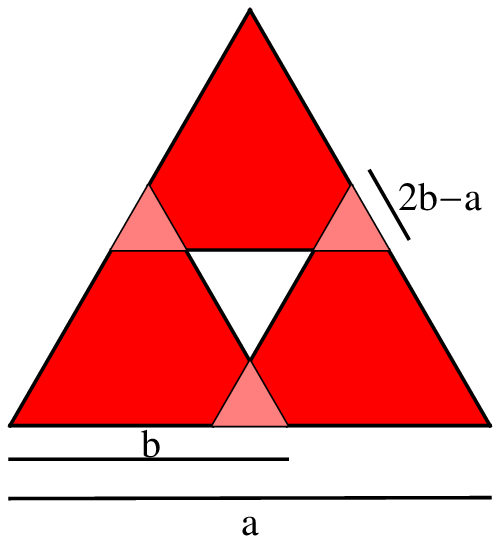}}
\caption{\label{fig:sqrt3} Geometric proof of the irrationality of $\sqrt{3}$}
\end{center}\end{figure}

We generalize Tennenbaum's geometric proof to show $\sqrt{3}$ is irrational. Suppose not, so $\sqrt{3} = a/b$, and again we may assume that $a$ and $b$ are the smallest integers satisfying this. This time we have $a^2 = 3b^2$, which we interpret as the area\footnote{Note that the area of an equilateral triangle is proportional to the square of its side $s$; specifically, the area is $s^2 \cdot \sqrt{3}/4$.} of one equilateral triangle of side length $a$ equals the area of three equilateral triangles of side length $b$. We represent this in Figure \ref{fig:sqrt3}, , which consists of three equilateral triangles of side length $b$ placed at the corners of an equilateral triangle of side length $a$. Note that the area of the three doubly covered, pink triangles (which have side length $2b-a$) is therefore equal to that of the uncovered, equilateral triangle in the middle. This is clearly a smaller solution, contradiction!\footnote{For completeness, the length of the equilateral triangle in the middle is $2a-3b$. To see this, note the dark red triangles have sides of length $b$, each pink triangle has sides of length $2b-a$ and thus the middle triangle has sides of length $b - 2(2b-a)=2a-3b$. This leads to $3(2b-a)^2 = (2a-3b)^2$, which after some algebra we see is a smaller solution.}

%%%%%%%%%%%%%%%%%%%%%%%%%%%%%%%%%%%%%%%%%%%%%%%%%%%%%%%%%%%%%%%%%%%%%%%%%%%%%%%%%%%%%%%%%%%%%%%%%%%%%%%%%%%%%%%%

%%%%%%%%%%%%%%%%%%%%%%%%%%%%%%%%%%%%%%%%%%%%%%%%%%%%%%%%%%%%%%%%%%%%%%%%%%%%%%%%%%%%%%%%%%%%%%%%%%%%%%%%%%%%%%%%

\section{The square-root of $5$ is irrational}

For the irrationality of $\sqrt{5}$, we have to slightly modify our approach as the overlapping regions are not so nicely shaped. As the proof is similar, we omit some of the details. Similar to the case of $\sqrt{3}$ and triangles, there are proportionality constants relating the area to the square of the side lengths of regular $n$-gons; however, as these constants appear on both sides of the equations, we may ignore them.

\begin{figure}[h]
\begin{center}
\scalebox{.7}{\includegraphics{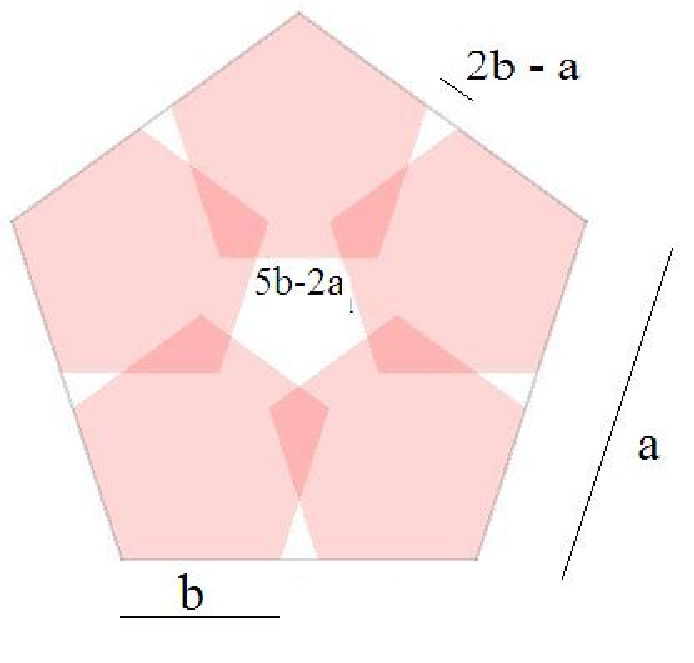}}
\caption{\label{fig:sqrt5Fig1} Geometric proof of the irrationality of $\sqrt{5}$. }
\end{center}\end{figure}

Suppose $a^2 = 5b^2$ with, as always, $a$ and $b$ minimal. We place five regular pentagons of side length $b$ at the corners of a regular pentagon of side length $a$ (see Figure \ref{fig:sqrt5Fig1}).  Note that this gives five small triangles on the edge of the larger pentagon which are uncovered, one uncovered regular pentagon in the middle of the larger pentagon, and five kite-shaped doubly covered regions. As before, the doubly covered region must have the same area as the uncovered region.

\begin{figure}
\begin{center}
\scalebox{.6}{\includegraphics{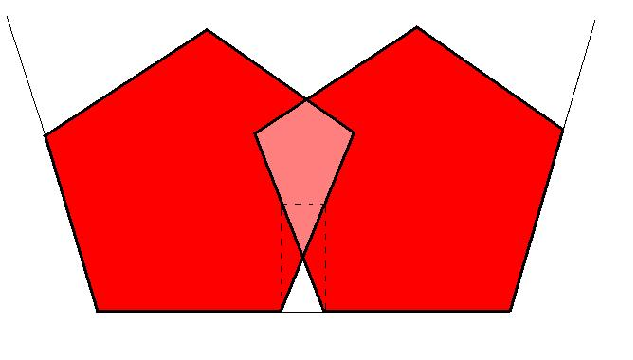}}\ \scalebox{.6}{\includegraphics{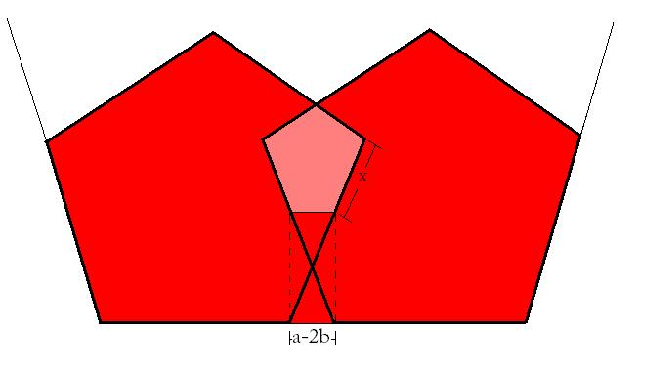}}
\caption{\label{fig:sqrt5Fig23} Geometric proof of the irrationality of $\sqrt{5}$: the kites, triangles and the small pentagons. }
\end{center}\end{figure}

We now take the uncovered triangles from the edge and match them with the doubly covered part at the ``bottom'' of the kite, and regard each as covered once instead of one covered twice and one uncovered (see Figure \ref{fig:sqrt5Fig23}).
This leaves five doubly covered pentagons, and one larger pentagon uncovered. We show that the five doubly covered pentagons are all regular, with side length $a-2b$. From this we get that the side lengths of the middle pentagon are all equal to $b - 2(a-2b) = 5b-2a$; later we'll show that the five angles of the middle pentagon are all equal, and thus it too is a regular pentagon. We now have a smaller solution, namely $5(a-2b)^2 = (5b-2a)^2$ (note that $a-2b < b$, as $a = b\sqrt{5}$, and $2 < \sqrt{5} < 3$), and thus we will have our contradiction.

To show the smaller pentagons are regular, we just need to show two things: (1) all the angles are $3\pi/5$; (2) the lower left and lower right sides (which have equal side length by symmetry, which we denote by $x$) in Figure \ref{fig:sqrt5Fig23} have length $a-2b$, which is the length of the bottom side. If we can show these three sides are all equal, then since all five angles are equal the remaining two sides must also have length $a-2b$ and thus the pentagon is regular.\footnote{This is similar to the ASA, or angle-side-angle, principle for when two triangles are equal. There is a unique pentagon once we specify all angles and three consecutive sides; as a regular pentagon satisfies our conditions, it must be the only solution.}

For (1), the sum of the angles of a pentagon\footnote{To see this, connect each corner to the center. This forms 5 triangles, each of which gives $\pi$ radians. We must subtract $2\pi$ for the sum of the angles around the center, which gives $3\pi$ for the sum of the pentagon's angles.} is $3\pi$, and for a regular pentagon each angle is $3\pi/5$. The two base angles of the triangles are thus $2\pi/5$ (as they are supplementary angles), and thus the two angles in the smaller pentagon next to the triangle's base are also $3\pi/5$. The two angles adjacent to these are just internal angles of a regular pentagon, and thus also $3\pi/5$. As the sum of all the angles is $3\pi$ and we've already accounted for $4 \cdot 3\pi/5$, this forces the top angle to be $3\pi/5$ as desired. The proof that the angles in the middle pentagon are all equal proceeds similarly (or by symmetry).

We now turn to (2). The length of the left and right sides of the originally uncovered triangle at the bottom of Figure \ref{fig:sqrt5Fig23} is $\frac{(a-2b)/2}{\cos(2\pi/5)}$; to see this, note the base angles of the triangle are $\pi - \frac{3\pi}{5}$ (as each angle of a regular pentagon is $3\pi/5$ and the triangle's base angles are supplementary), and by standard trigonometry ($\cos({\rm angle}) = {\rm adjacent}/{\rm hypothenuse}$) we have $x \cos(\pi - \frac{3\pi}{5}) = \frac{a-2b}{2}$. Using $a/b = \sqrt{5}$, we see $x = b - \frac{2(a-2b)}{2\cos(2\pi/5)} = b(1-\frac{\sqrt{5}-2}{\cos(2\pi/5)})$. Noting the formula $\cos(2\pi/5) = \frac{1}{4}(-1+\sqrt{5})$, the previous expression for $x$ simplifies to $x = b(\sqrt{5}-2) = a-2b$ (as we are assuming $\sqrt{5}=a/b$), and the argument is done.

%%%%%%%%%%%%%%%%%%%%%%%%%%%%%%%%%%%%%%%%%%%%%%%%%%%%%%%%%%%%%%%%%%%%%%%%%%%%%%%%%%%%%%%%%%%%%%%%%%%%%%%%%%%%%%%%

%%%%%%%%%%%%%%%%%%%%%%%%%%%%%%%%%%%%%%%%%%%%%%%%%%%%%%%%%%%%%%%%%%%%%%%%%%%%%%%%%%%%%%%%%%%%%%%%%%%%%%%%%%%%%%%%

\section{How far can we generalize: Irrationality of $\sqrt{6}$}

\begin{figure}
\begin{center}
%\scalebox{.6}{\includegraphics{IrrationalitySqrt6.eps}}
\scalebox{.6}{\includegraphics{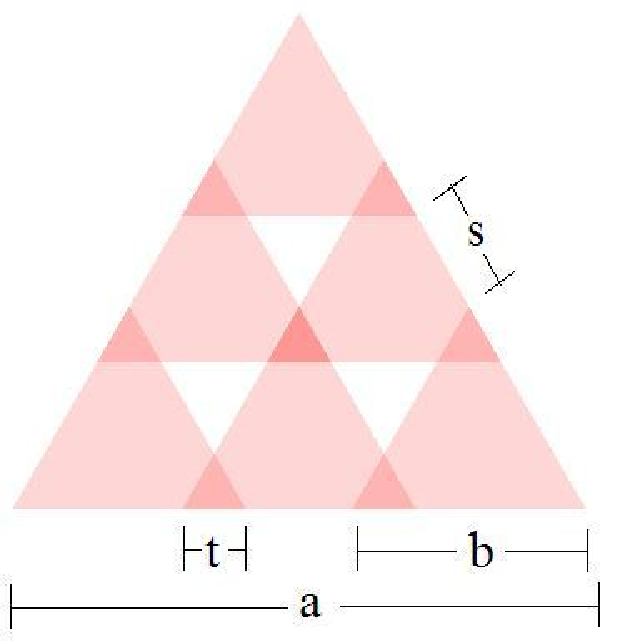}}
\caption{\label{fig:sqrt6} Geometric proof of the irrationality of $\sqrt{6}$.}
\end{center}\end{figure}

We conclude with a discussion of one generalization of our method that allows us to consider triangular numbers.
Figure \ref{fig:sqrt6}
 shows the construction for the irrationality of $\sqrt{6}$. Assume $\sqrt{6} = a/b$ so $a^2=6b^2$. The large equilateral triangle has side length $a$ and the six medium equilateral triangles have side length $b$. The 7 smallest equilateral triangles (6 double counted, one in the center triple counted) have side length $t = (3b-a)/2$. It's a little work, but not too bad, to show the triple counted one is the same size. For the three omitted triangles, they are all equilateral (angles equal) and of side length $s = b - 2(3b-a)/2 = a-2b$. As the area of the smaller equilateral triangles is proportional to $t^2$ and for the larger it is proportional to $s^2$, we find $8t^2 = 3s^2$ or $16t^2 = 6s^2$ so $(4t/s)^2 = 6$. Note that although $t$ itself may not be an integer, $4t = 2(3b-a)$ is an integer. We claim $4t < a$ and $s < b$, so that this is a smaller solution. Clearly $s<b$, and as the ratio is $\sqrt{6}$, the other claim now follows.

Can we continue this argument? We may interpret the argument here as adding three more triangles to the argument for the irrationality of $\sqrt{3}$; thus the next step would be adding four more triangles to this to prove the irrationality of $\sqrt{10}$. Proceeding along these lines lead us to study the square-roots of triangular numbers.\footnote{Triangular numbers are of the form $n(n+1)/2$ for some positive integer $n$. Note the first few are $1, 3, 6, 10, 15, \dots$.} We continue more generally by producing images like Figure \ref{fig:sqrt6} with $n$ equally spaced rows of side length $b$ triangles. This causes us to start with $a^2 = \frac{n(n+1)}{2}\ b^2$, so we can attempt to show that $\sqrt{n(n+1)/2}$ is irrational.

By similar reasoning to the above, we see that the smaller multiply covered equilateral triangles all have the same side length $t$, and that the uncovered triangles also all have the same side length $s$. Further $t$ equals $(nb-a)/(n-1)$, and we have that $s = b-2t$, so $s = b - 2(nb-a)/(n-1) = (2a-(n+1)b)/(n-1)$. To count the number of side length $t$ triangles, we note that there will be $(n-2)(n-1)/2$ triply covered triangles (as there is a triangle-shaped configuration of them with $n-2$ rows), and that there will be $3(n-1)$ doubly covered triangles around the edge of the figure, for a grand total of $2(n-2)(n-1)/2 + 3(n-1) = (n-1)(n+1)$ coverings of the smaller triangle. Further, note that in general there will be $(n-1)n/2$ smaller, uncovered triangles, so we have that $(n-1)(n+1)t^2 = ((n-1)n/2)s^2$. Writing out the formula for $s, t$ (to verify that our final smaller solution is integral), we have $(n-1)(n+1)((nb-a)/(n-1))^2 = ((n-1)n/2)((2a-(n+1)b)/(n-1))^2$. We now multiply both sides of the equation by $n-1$ to ensure integrality, giving $(n+1)(nb-a)^2 = (n/2)(2a-(n+1)b)^2$. We multiply both sides by $n/2$ to achieve a smaller solution to $a^2 = (n(n+1)/2) b^2$, giving us $(n(n+1)/2)(nb-a)^2 = (n(2a-(n+1)b)/2)^2$. Note that this solution is integral, as $n$ odd implies that $2a - (n+1)b$ is even. Finally, to show that this solution is smaller, we just need that $nb-a < b$. This is equivalent to $n-\sqrt{n(n+1)/2} < 1$.

We see that this inequality holds for $n \leq 4$, but not for $n > 4$. So, we have shown that the method used above to prove that $\sqrt{6}$ is irrational can also be used to show that $\sqrt{10}$ (the square root of the fourth triangular number) is irrational, but that this method will not work for any further triangular numbers. It is good (perhaps it is better to say, `it is not unexpected') to have such a problem, as some triangular numbers are perfect squares. For example, when $n=49$ then we have $49\cdot 50/2 = 7^2 \cdot 5^2$, and thus we should not be able to prove that this has an irrational square-root!

%%%%%%%%%%%%%%%%%%%%%%%%%%%%%%%%%%%%%%%%%%%%%%%%%%%%%%%%%%%%%%%%%%%%%%%%%%%

\ \\

\end{document}